
\input amstex
\documentstyle{amsppt}
\pagewidth{6.5truein}
\pageheight{9truein}
\overfullrule=0in
\def\inv#1{\roman{inv}(#1)}
\def\evalat{\bgroup\left|\,\smallmatrix\format\l&\quad\l\\}
\def\endevalat{\endmatrix\right.\egroup}

\def\qqline#1#2{\evalat {}\cr {}\cr#1\hfill\cr #2 \hfill\cr\endevalat }
\def\qqqline#1#2#3{\evalat {}\cr{}\cr #1\hfill\cr
 #2 \hfill\cr #3 \hfill\cr\endevalat}
\def\qline{\qqqline{q_k{\rightarrow }qq_k,\ 2\leq k{<}i}{q_k{\rightarrow }q_{k{+}1},\ i\leq
k{<}n}{q{\rightarrow }q}}
\def\q1line#1{\qqqline{q_k{\rightarrow}qq_k,\ 2\leq k{<}#1}{q_k{\rightarrow}q_{k{+}1},\ #1\leq
k{<}n}{q{\rightarrow}q}}
\def\pdir#1{{\partial\over{\partial #1}}}
\def\p2dir#1{{{\partial^2}\over{\partial #1^2}}}

\topmatter
\title The Enumeration of Permutations With a Prescribed Number of 
``Forbidden'' Patterns
\endtitle
\author John Noonan and Doron Zeilberger
\endauthor
\address  Department of Mathematics, Temple University, Philadelphia, PA
       19122\endaddress
\email noonan{\@}math.temple.edu, zeilberg{\@}math.temple.edu, 
{\it WWW homepage}: http://www.math.temple.edu/\~{}noonan, 
http://www.math.temple.edu/\~{}zeilberg\endemail
\abstract
We initiate a general approach for the fast enumeration of permutations with
a prescribed number of occurrences of `forbidden' patterns, that
seems to indicate that the enumerating sequence is always P-recursive. 
We illustrate the
method completely in terms of the patterns `abc',`cab' and `abcd'.

\endabstract
\endtopmatter
\document
\baselineskip=17pt

\head 0. Introduction\endhead

The only increasing permutation on $\{1,2, \dots ,n \}$, $[1,2,\dots,n]$, 
has the
property that it has no decreasing subsequence of length $2$,
i.e. there are no $i$ and $j$ such that $1 \leq i < j \leq n$, and
$\pi[i]>\pi[j]$. Thus one measure of how scrambled a permutation is,
is its {\it number of inversions}, $\inv \pi$, which is the number
of occurrences of the `pattern' $ba$.

This can be generalized to an arbitrary set of patterns.
Given a permutation on $\{1, \dots , n\}$, 
we define a {\it pattern} as 
a permutation on $\{1, \dots ,r \}$ where $r\leq n$. For $r\leq n$, we say that a permutation 
$\sigma\in S_n$ has the pattern $\pi\in S_r$ if there exist 
$1\leq i_1< i_2<\dots<i_r\leq n$ such that 
$\pi\equiv [\sigma(i_1),\sigma(i_2),\dots,\sigma(i_r)]$ in reduced form. 

The {\it reduced form} of a permutation $\sigma$ on a set  
$\{j_1,j_2,\dots,j_n\}$ where $j_1<j_2<\dots<j_n$ is the permutation 
$\sigma_1\in S_n$ obtained by renaming the objects of the permutation $\sigma$ in 
the obvious way, so that $j_1$ is renamed 1 and $j_2$ is renamed 2 and so on.  
Thus the reduced form of the permutation $25734$ is $14523$ and the reduced 
form of $579$ is $123$.

To simplify things, when discussing a particular pattern, we will use 
its alphabetic equivalent. Thus an increasing subsequence of length 3 
is an $abc$ pattern, an inversion is a $ba$ pattern and a 
decreasing subsequence of length 4 is a $dcba$ pattern. Other 
patterns we will discuss include $abcd$, $bac$ and $cab$.\par
The number of permutations which contain no increasing subsequence of 
length three (i.e. $abc$ avoiding) is known to be 
$C_n:={1\over {n+1}}\binom{2n}{n}$, the Catalan 
numbers. It is also known \cite{6} 
that given any pattern of length three, the number 
of permutations avoiding that pattern is also $C_n$. 

Herb Wilf raised the question:
For any pattern $\pi$, what can you say about $a_\pi(n)$, the number of 
permutations on $\{1\dots n\}$ that avoid the pattern $\pi$?

It follows from the Robinson-Schenstead algorithm 
and the hook-length formula\cite{3} that for any $r$,
the number of permutations with no increasing subsequence of length $r$,
 is a certain binomial-coefficient multisum, from which
it follows immediately \cite{8} that it is P-recursive
(holonomic) (i.e., it satisfies a linear recurrence with polynomial
coefficients in $n$).

A natural conjecture is: For any given finite set of patterns, $PAT$,  the sequence
$$
a_{PAT}(n):= \biggl|\Bigl\{\sigma\in S_n :\ \sigma 
\ \roman{has\ no\ occurrences\ of\ the\ given\ patterns}\Bigr\}\biggr|
$$
is $P$-recursive. More generally, for any such set of patterns 
$PAT = \{ pat_1, pat_2, \dots, pat_r\}$
and any specified sequence of integers $m_p$, one for each pattern $p \in PAT$,
$$
a_{PAT}^{\{m_1,m_2,\dots,m_r\}}(n):= \biggl|\Bigl\{\sigma\in S_n :\ \sigma \roman{\ has\ \italic{exactly}}
\ m_i\ \roman{ occurrences\ of\ the\ pattern\ }pat_i,
\ 1\leq i\leq r \Bigr\}\biggr|,
$$
\penalty10000
is P-recursive.

In this paper, we will present a method for the `fast' (polynomial
time in $n$) enumeration of such sequences $a(n)$, that seems 
to support our conjecture
that it should always be P-recursive (holonomic) in $n$.

Until later in the paper, we will consider only one pattern, as the main ideas
are already present there. 

The natural object is the {\it generating function:}
\penalty10000$$
F_n^\pi(q)=\sum_{\sigma\in S_n}q^{\varphi(\sigma)}=\sum c_i(n) q^i
$$\penalty10000
where $\varphi(\sigma)$ denotes the number of 
subsequences, $x_1x_2\dots x_r$, present in $\sigma$ that 
reduce to the given pattern $\pi$. 
We will later show how to derive
a recursive functional equation, from which it should be possible to
extract efficient recurrences for the coefficients $c_i(n)$ of $F_n^{\pi}(q)$ 
for small $i$.
As we said above, our approach seems to indicate that each $c_i(n)$ is
$P$-recursive in $n$. It is easy to see that it cannot be
also $P-$recursive in $i$.

Another possibility is to expand the polynomials around $q=1$, so 
that we have 
$$
F_n^\pi(q)=\sum_{i=0}^{degree}b_i(n)(q-1)^i.
$$
Here $b_0(n)$ is the total number of permutations, $n!$. $b_1(n)$ is the total 
number  of
$\pi=x_1x_2\dots x_r$ patterns present in all the
permutations of set $S_n$. The average number of 
$x_1x_2\dots x_r$ patterns present in the permutations on $\{1\dots n\}$ would
simply be $\displaystyle {{b_1(n)}\over {n!}}$. It is easy to see \cite{9} that
it is always a polynomial in $n$, as are all the other coefficients
$b_i(n)$, for each fixed $i$.
Thus, we can compute the first few coefficients  $b_i(n)$ of $F_n^{\pi}(q+1)$
by brute force. However to get the full $F_n^{\pi}(q)$ we would need
the full $F_n^{\pi}(q+1)$, and the coefficients $b_i(n)$, while always
polynomials in $n$, get increasingly complicated as $i$ grows bigger.

Instead we will answer a more general 
question: How many permutations on $\{1\dots n\}$ avoid $\pi=x_1x_2\dots x_r$ 
{\it and} $\{ \pi_i\}_{i=1}^m$ where $\{\pi_i\}$ is a set of other forbidden
patterns for which one or more of the entries in the forbidden 
subsequence is specified. Among the $\pi_i$ might be the pattern $ab4$ 
which simply describes an $abc$ pattern in which the last entry is $4$.
Another might be a $3b$ pattern which would be a non-inversion in which 
the first entry is $3$.\par 

Once we have such a method, we would also like to compute the number of 
permutations of a given size which avoid a given set of patterns and 
compute the number of permutations of a given size which have occurences of
patterns from the given
set, a prescribed number of times. We may ask how many 
permutations on $\{1\dots n\}$ avoid both $abc$ and $cab$?; how many 
avoid $abc$ and have exactly 1 $cab$?; etc.

\head 1 Counting 
Permutations with a Presribed Number of $abc$ Patterns \endhead
\head 1.0 Definitions\endhead
\definition{Definition 1.1} Given $\sigma\in S_n$, an $abc$ pattern is a
sequence $i,j,k$ where $0<i<j<k\leq n$ and $\sigma(i)<\sigma(j)<\sigma(k)$.
\enddefinition
\definition
{Definition 1.2}
 $\varphi_{abc}(\sigma)$:= the number of $abc$ patterns of $\sigma$.
\enddefinition

For example $\varphi_{abc}(4321)=0$, $\varphi_{abc}(1234)=4$,
and $\varphi_{abc}(2314)=1$.

\definition{Definition 1.3} Given $\sigma\in S_n$, an $aj$ pattern is a 
sequence $i,k$ where $0<i<k\leq n$ and $\sigma(i)<\sigma(k)=j$.\enddefinition

For example $\varphi_4 (15342)=2$.

\definition{Definition 1.4} $\varphi_j(\sigma) := $ the number of $aj$
patterns of $\sigma$.
\enddefinition
\definition{Definition 1.5} $P^{(r)}(n,I)$ is the number of permutations on 
$\{1\dots n\}$ with exactly $r$ $abc$ patterns and no $aj$ patterns for 
$j\leq I$.\enddefinition
We will use $P(n,I)$ to denote $P^{(0)}(n,I)$, i.e. the number of permutations on $\{1\dots n\}$ with 
no $abc$ patterns and no $aj$ patterns for $j\leq I$.\par
Using the above definitions, the polynomial $F_n^{\pi}(q)$ described earlier would for 
this example be defined as 
$$
F_n^{123}(q):=\sum_{\sigma\in S_n}q^{\varphi_{abc}(\sigma)}.
$$

So, e.g., $F_1^{123}(q)=1$ ,\quad $F_2^{123}(q)=2$ ,\quad $F_3^{123}(q)=5+q$ , \quad $F_4^{123}(q)=14+6q+3q^2+q^4$.

\definition{Definition 1.6} For $\sigma\in S_n$, define:
$$
wt(\sigma) :=q^{\varphi_{abc}(\sigma)} q_2^{\varphi_2(\sigma)}q_3^{\varphi_3(\sigma)}q_4^{\varphi_4(\sigma)}\dots 
q_n^{\varphi_n(\sigma)}.
$$
\enddefinition

For example, $\displaystyle wt(2314)=q^1 q_2^0 q_3^1 q_4^3=qq_3q_4^3$.

\definition{Definition 1.7} 
$$
P_n(q_2,q_3,\dots,q_n;q):=\sum_{\sigma\in S_n}wt(\sigma).
$$
\enddefinition
$P_n$ is the `generalized' form of $F_n^{123}$. Comparing the two we see that $F_n^{123}(q)=P_n(1,1\dots,1;q)$.
\head 1.1 The functional equation\endhead
In order to illustrate the present method,
we will treat the simplest non-trivial case, by rederiving the 
well-known formula for the  
number of permutations on $\{1\dots n\}$ with no $abc$ 
patterns by first obtaining a method to compute the number of 
permutations on $\{1\dots n\}$ with no $abc$ patterns and no $aj$ patterns for $j<I$ where $I$ is any integer between 0 and $n$. We then merely set $I=0$ and obtain the desired sequence.\par

In order to find explicit or recursive descriptions for the coefficients of 
$F_n^{123}$, we will establish a recursive functional equation for $P_n$. 
Let $\sigma\in S_n$, $\sigma(n)=i\neq 1$. Let $\sigma_1$ be the permutation on $\{1..i-1,i+1..n\}$ obtained by removing the last entry of $\sigma$, i.e.
$$
\sigma_1(j):=\sigma(j), \ 1\leq j \leq n-1.
$$
 Then 
$$
wt(\sigma)=wt(\sigma_1)q_i^{i-1}q^{\sum_{j=2}^{i-
1}\varphi_j(\sigma_1)},
$$
hence
$$
wt(\sigma)=q^{\varphi_{abc}(\sigma)}\prod_{j=1}^nq_j^{\varphi_j(\sigma)}=q^{abc(\sigma_1)}q_i^{i-1}\prod_{j<
i}(qq_j)^{\varphi_j(\sigma_1)}\prod_{j>i}q_j^{\varphi_j(\sigma_1)}.\eqno(1)
$$

If $\sigma(n)=1$, then defining $\sigma_1$ as above with $i=1$, we have
$$
wt(\sigma)=wt(\sigma_1).
$$ 
Summing over all $\sigma\in S_n$ on the left hand side of (1) is equivalent 
to summing first over $i$, then over $\sigma_1\in S_{n-1}$ on the right. Making the necessary shift in variables to account for the fact that 
$\sigma_1$ above is a permutation on $[1,..,i-1,i+1,..,n]$, we 
have 
$$\eqalign{
\sum_{\sigma\in S_n}wt(\sigma)=\sum_{i=1}^n
\sum_{{\sigma(n)=i}\atop{\sigma\in S_n}}wt(\sigma)&=\sum_{\sigma_1\in S_{n{-}1}^{(1)}}wt(\sigma_1)\cr
&+\sum_{i=2}^n\sum_{\sigma_1\in S_{n-1}^{(i)}}q_i^{i{-}1}\left[\sum_{\sigma_1\in S_{n{-}1}}q^{abc(\sigma_1)}\prod_{j<
i}(qq_j)^{\varphi_j(\sigma_1)}\prod_{j>i}q_j^{\varphi_j(\sigma_1)}\right]\cr}
$$
where $S_{n-1}^{(i)}$ is the set of permutations on
$\{1,2,\dots,i-1,i+1,\dots,n-1,n\}$.

Hence
$$
P_n(q_2,q_3,\dots,q_n;q)= P_{n-1}(q_3,q_4,\dots,q_n;q)+\sum_{i=2}^n 
q_i^{i-1}P_{n-1}(qq_2,qq_3,\dots,qq_{i-1},q_{i+1},\dots,q_n;q).\eqno(2)
$$
This recurrence for $P_n$ is the basis for all that follows.

\head 1.2 A recurrence for $P(n,I)$\endhead

We note that from our definition of $P_n$, 
$P_{n}(\overbrace{0,0,\dots ,0}^{I-1},\overbrace{1,\dots,1}^{n-I-1};0)$ is 
the number of permutations on $\{1\dots n\}$ with no $abc$ patterns and no 
$aj$ patterns for $j\leq I$. So $P_n(1,1,\dots,1;0)$ is $a_{123}(n)$, the 
number 
of permutations on $n$ elements with no $abc$ patterns, and  
$P_{n}(\overbrace{0,0,\dots ,0}^{I-1},\overbrace{1,\dots,1}^{n-I-1};0)=P(n,I)$.
We first tackle the question of a recurrence for 
$a_{123}(n)=P_n(1,1,\dots,1;0)$. 
Using $(2)$, we have
$$
P_n(1,1,\dots,1;0)= P_{n-1}(1,1,\dots,1;0)+\sum_{i=2}^n 
P_{n-1}(\overbrace{0,0,\dots 
,0}^{i-2},\overbrace{1,\dots,1}^{n-i-1};0).\eqno(3)
$$
It would be too much to hope for a nice clean recurrence on the first 
try, we seem to have picked up a few uninvited guests, namely the 
$P_{n-1}(0,0,\dots ,0,{1,\dots,1};0)$. By $(2)$ again,
$$
P_{n}(\overbrace{0,0,\dots ,0}^{I-1},\overbrace{1,\dots,1}^{n-I-1};0) =
P_{n-1}(\overbrace{0,0,\dots 
,0}^{I-2},\overbrace{1,\dots,1}^{n-I-1};0)+\sum_{i=I+1}^n 
P_{n-1}(\overbrace{0,0,\dots 
,0}^{i-2},\overbrace{1,\dots,1}^{n-i-1};0).\eqno(4)
$$

In terms of $P(n,I)$, $(4)$ becomes
$$
P(n,I) = P(n-1,I-1) + \sum_{i=I+1}^n P(n-1,i-1) = \sum_{i=I}^n 
P(n-1,i-1).\eqno(5)
$$

The number of permutations on $n$ elements with no $abc$ patterns is 
$P_n(1,1,\dots,1;0)=P(n,1)=a_{123}(n)$. Unfortunately when we try to use 
$(5)$ to find $P(n,1)$ we are required to define $P(n,0)$. To do 
so, one needs only look as far as $(2)$. By $(2)$, we have
$$
P(n,1) = P(n-1,1) + \sum_{i=2}^n P(n-1,i-1).
$$

Comparing this with $(5)$ we see that we should 
define $P(n,0)=P(n,1)$. Using this definition, $(5)$ is valid 
for $n\geq 1, i> 0$. To complete the scheme, we need some form of 
initial conditions. These are readily supplied by the observation that 
the 
number of permutations on $n$ elements with no $abc$ patterns and no 
$aj$ patterns $2\leq j\leq n$ is 1, namely the permutation 
$[n,n-1,\dots, 3,2,1]$, so $P(n,n)= 1$.
\smallskip
Finally, we may simplify $(5)$ by examining\par
\noindent 
$P(n,I)-P(n,I+1)=\sum_{i=I}^n P(n-1,i-1)-\sum_{i=I+1}^n P(n-1,i-1)= 
P(n-1,I-1)$ which gives us 
$$
P(n,I) = \cases 1,\hfill &\roman{if}\ n=I\hfill\cr
P(n,1),\hfill &\roman{if}\ I=0\hfill\cr
P(n,I+1)+P(n-1,I-1),&\roman{otherwise}\cr \endcases \eqno(6)
$$
Using this recurrence, we may quickly generate a large number of 
$P(n,I)$.
\bigskip
\centerline{
\vbox{
\halign{#\hfil\qquad& \qquad\hfil#\qquad& \hfil#\qquad& \hfil#\qquad& \hfil#\qquad& \hfil#\qquad& \hfil#\qquad& \hfil#\qquad& \hfil#\qquad\cr
\multispan{9}\hfil\bf{Table 1}\hfil\cr
\multispan{9}\hfil Values of $P(n,I)$\hfil\cr
\noalign{\smallskip\hrule\smallskip}
n&I=0&1&2&3&4&5&6&7\cr
\noalign{\smallskip\hrule\smallskip}
0& 1& & & & & & & \cr
1& 1& 1& & & & & & \cr
2& 2& 2& 1& & & & &\cr
3& 5& 5& 3& 1& & & &\cr
4& 14& 14& 9& 4& 1& & &\cr
5& 42& 42& 28& 14& 5& 1& &\cr
6& 132& 132& 90& 48& 20& 6& 1& \cr
7& 429& 429& 297& 165& 75& 27& 7& 1\cr
8& 1430& 1430& 1001& 572& 275& 110& 35& 8\cr
9& 4862& 4862& 3432& 2002& 1001& 429& 154& 44\cr
10& 16796& 16796& 11934& 7072& 3640& 1638& 637& 208\cr
\noalign{\smallskip\hrule}}}}
\bigskip

This enables us to conjecture, and immediately prove
(by verifying $(6)$ and the initial conditions), the closed form 
$P(n,I)= {{I+1}\over{n+1}}\binom{2n-I}{n}$, the celebrated 
ballot numbers\cite{3}. Evaluating at $I=1$ yields yet another proof of
the well known fact that the number of permutations on $\{1\dots n\}$ with 
no $abc$ patterns, $a_{123}(n)$, equals $C_n$, the Catalan number. 
This proof is longer and far less elegant than the combinatorial proofs
of \cite{6}. Its only virtue is that it illustrates a {\it general} method.
\medskip

\head 1.3 The number of permutations with exactly one $abc$ pattern.\endhead

Recently one of us \cite{4} proved that 
the number of permutations on $\{1\dots n\}$ with exactly one $abc$ 
pattern is ${3\over n}\binom{2n}{n+3}$. We will now present an alternative
proof using the present method.

The number of permutations with exactly one $abc$ pattern is the 
constant term of the derivative of the 
function $F=\sum_{\sigma\in S_n}q^{\varphi_{abc}(\sigma)}$. 
By differentiating $(2)$  with respect to $q$, we find a recurrence for 
$\displaystyle{{{\partial 
}\over{\partial q}}P_{n}}$,
$$
\pdir qP_n(q_2,\dots,q_n;q)= 
\pdir qP_{n-1}(q_3,\dots,q_n;q)+\sum_{i=2}^n 
q_i^{i-1}{\partial\over{\partial 
q}}P_{n-1}(qq_2,qq_3,\dots,qq_{i-1},q_{i+1},\dots,q_n;q). \eqno(7)
$$
Here things begin to get a bit tricky and we will need to employ the 
chain rule. We will evaluate the right side 
of the above equation in terms of partial derivatives with respect to
the 
positions that $q$ occupies. For example, 
$$
\pdir qP_3(qq_2,q_3;q) = 
q_2{\partial\over{\partial 
q_2}}P_3\qqqline{q_2\rightarrow qq_2}{q_3\rightarrow q_3}{q\rightarrow q}+{\partial\over{\partial 
q}}P_3\qqqline{q_2\rightarrow qq_2}{q_3\rightarrow q_3}{q\rightarrow q}.$$

The preceding 
notation means: ``first find the partial derivative of the function $P_3$
then make the necessary substitutions for $q_i$." Continuing we have
$$\eqalignno{
\quad\pdir qP_n(q_2,&q_3,\dots,q_n;q)= 
{\partial\over{\partial
q}}P_{n{-}1}(q_3,q_4,\dots,q_n;q)&\cr&+\sum_{i=2}^n q_i^{i{-
}1}\left[ \pdir qP_{n{-}
1}\qline{+}\sum_{j=2}^{i{-}1}q_j{\partial\over{\partial 
q_j}}P_{n{-}1}\qline\right].\quad&(8)}
$$
We are now ready to tackle  $\displaystyle{{{\partial 
}\over{\partial q}}P_{n}\qqline{q\rightarrow 0}{q_i\rightarrow 1,\  2{\leq} i{\leq} n}}$

$$\eqalign{
{{\partial 
}\over{\partial q}}P_{n}&\qqline{q_k{\rightarrow }1,\  2{\leq} k{\leq} n}{q{\rightarrow }0}\cr
&= {{\partial 
}\over{\partial q}}P_{n{-}1}\qqline{q_k{\rightarrow }1,\  2{\leq} k{\leq}
n{-}1}{q{\rightarrow }0} {+} 
\sum_{i{=}2}^n \left[
{{\partial 
}\over{\partial q}}P_{n{-}1}\qqqline{q_k{\rightarrow }0,\  2{\leq}
k{<}i}{q_k{\rightarrow }1,\ i{\leq} 
k{\leq} n{-}1}{q{\rightarrow }0}
{+}\sum_{j{=}2}^{i{-}1}{{\partial 
P_{n{-}1}}\over{\partial q_j}}\qqqline{q_k{\rightarrow }0,\ 2{\leq}
k{<}i}{q_k{\rightarrow }1,\ i{\leq} 
k{\leq} n{-}1}{q{\rightarrow }0}
\right].\cr}
$$
It is readily seen that $P^{(1)}(n,I)$,
the number of permutations on $\{1\dots n\}$ with exactly one 
$abc$ pattern and  no $aj$ patterns for $j\leq I$ can be expressed as:
$$
P^{(1)}(n,I) ={{\partial 
}\over{\partial q}}P_{n}\qqqline{q_k\rightarrow 0,\  2\leq k<I+1}{q_k\rightarrow 1,\ I+1\leq
k\leq
n}{q\rightarrow 0}.
$$
Furthermore, to simplify notation, let
$$
P^{(1)}_j(n,I):=\displaystyle{{\partial 
P_n}\over{\partial q_j}}\qqqline{q_k\rightarrow 0,\  2\leq k<I+1}{q_k\rightarrow 1,\ I+1\leq
k\leq
n}{q\rightarrow 0}.
$$
The actual meaning of $P^{(1)}_j(n,I)$ is immaterial,
but the astute reader will see that combinatorially, $P^{(1)}_j (n,I)$ is 
the number of permutations on $\{1\dots n\}$ with no $abc$ patterns, no 
$ak$ patterns for $k\leq I,\ k\neq j$, and exactly 1 $aj$ 
pattern if $j\leq I$ or {\it at least} one $aj$ pattern if 
$j>I$. The 
recurrence for
$P^{(1)}(n,I)$ follows from $(8)$:
$$P^{(1)}(n,I)=P^{(1)}(n-1,I-1)+\sum_{i=I+1}^n
\left[P^{(1)}(n-1,i-1)+\sum_{j=I+1}^{i-
1}P_j^{(1)}(n-1,i-1)\right].\eqno(9)
$$
\bigskip
Disregarding for a moment that we do not yet know what $P_j^{(1)}(n,I)$
is, we can state `initial' conditions for this recurrence. First we should
define $P^{(1)}(n,0) = P^{(1)}(n,1)$.  We can easily compute $P^{(1)}(n,n-2)$. $P^{(1)}(n,n-2)$ is the number of permutations on $\{1\dots n\}$ with exactly 1 $abc$ pattern and no $aj$ 
patterns for $j\leq n-2$. If $\sigma$ is one such permutation then $\sigma$ 
has the form $[n-2,n-1,n-3,\dots,n-i,n,n-i-1,\dots,2,1]$, for some $i$. 
There are exactly $n-2$ such permutations, hence $P^{(1)}(n,n-2)=n-2$. Our `initial' conditions (perhaps they should be called
`boundary' conditions) for this recurrence are
$P^{(1)}(n,0)=P^{(1)}(n,1),\ P^{(1)}(n,n-2)=n-2$. \par
\smallskip

To obtain a recurrence for $P^{(1)}_j(n,I)$ we must return to $(2)$ 
and take the
partial derivative with respect to $q_j$. We begin for the case when 
$j\geq3$. 
$$
\eqalign{
{{\partial}\over{\partial q_j}}P_n(&q_2,\dots,q_n;q)\cr
&= {{\partial}\over{\partial q_{j{-}1}}}P_{n-
1}\qqline{q_k\rightarrow q_{k+1}}{q\rightarrow q}+\sum_{i=2}^{j{-}1} 
q_i^{i-1}{{\partial}\over{\partial q_{j{-}1}}}P_{n-1}\qqqline{q_k\rightarrow qq_k,\
2\leq
k<i}{q_k\rightarrow q_{k+1},\ i<k\leq n-1}{q\rightarrow q}\cr
&{+}\sum_{i=j+1}^{n} 
qq_i^{i{-}1}{{\partial}\over{\partial q_{j}}}P_{n-1}\qqqline{q_k\rightarrow qq_k,\
2\leq
k<i}{q_k\rightarrow q_{k{+}1},\ i<k\leq n{-}1}{q\rightarrow q}{+}(j{-}1)q_j^{j{-}2}P_{n{-
}1}(qq_2,\dots,qq_{j{-}
1},q_{j{+}1},\dots,q_n;q)\ .\cr}
$$

For $j\geq 3$ it follows that 
$$
P_j^{(1)}(n,I)=P_{j{-}1}^{(1)}(n-1,I-1)+\sum_{i=I+1}^{j{-}1}P_{j{-}1}^{(1)}(n-1,i-1)
+(j{-}1)\chi(j>I)P(n-1,j-1).
$$
\noindent
where $\chi(statement)=\left\{ {{\scriptstyle 0,\ statement\ \roman{is\ 
false}}\atop {\scriptstyle 1,\
\ statement\ \roman{is\ true}}}\right.$. We see that in $(9)$, we only need to know 
the values of
$P_j^{(1)}(n,I)$ for which $j\leq I$ so we have 
$$
P_j^{(1)}(n,I)=P_{(j-1)}^{(1)}(n-1,I-1).\eqno(10)
$$
When $j=2$ we have 
$$
\eqalignno{P_2^{(1)}(n,I)&= 
{{\partial 
}\over{\partial q_2}}P_{n}\qqqline{q_k\rightarrow 0,\  2\leq k< I+1}{q_k\rightarrow 1,\ I\leq
k\leq n}{q\rightarrow 0}
= P_{n-1}(\overbrace{0,\dots,0}^{I-2},\overbrace{1,\dots,1}^{n-I};0)&\cr
&=P(n-1,I-1)&(11)\cr}
$$
 Using (11), the recurrence
 $(10)$ simplifies to 
$$
\eqalign{
P_j^{(1)}(n,I)&=P_{j{-}1}^{(1)}(n-1,I-1)=P_{j{-}2}^{(1)}(n-2,I-2)
=\dots=P_2^{(1)}(n-j+2,I-j+2)\cr
&=P(n-j+1,I-j+1).}
$$

Now we can summarize our results so far in a recurrence for
$P^{(1)}(n,I)$ (see $(9)$).
$$
P^{(1)}(n,I)= \sum_{i=I}^{n}\left[ 
P^{(1)}(n-1,i-1)+\sum_{j=I+1}^{i-1}P(n-j,i-
j)\right ].\eqno(12)
$$

From $(12)$ it follows that
$$
\eqalign{
P^{(1)}(n,I)-P^{(1)}(n,I+1)&=P^{(1)}(n-1,I-1)+\sum_{i=I}^n \left(\sum_{j=I+1}^{i-1}
P(n-j,i-j)\right.- \left. \sum_{j=I+2}^{i-1}P(n-j,i-j)\right)\cr
&=P^{(1)}(n-1,I-1)+\sum_{i=I+2}^n P(n-I-1,i-I-1)\cr
&=P^{(1)}(n-1,I-1)+P(n-I,2).\cr}
$$

The last equation follows from $(5)$. Finally we have the 
following recursion for
$P^{(1)}(n,I)$:
$$
P^{(1)}(n,I)=\cases P^{(1)}(n,1),&\roman{if}\  I=0\cr n-2,&\roman{if}\ I=n-2\cr P^{(1)}(n,I+1)+P^{(1)}(n-1,I-1)+P(n-I,2),& \roman{otherwise}\endcases
$$
\bigskip
\centerline{
\vbox{
\halign{#\hfil\qquad& \qquad\hfil#\qquad&\hfil#\qquad& \hfil#\qquad& 
\hfil#\qquad& \hfil
#\qquad& \hfil#\qquad& \hfil#\qquad& \hfil#\qquad& \hfil#\qquad& \hfil#\qquad& 
\hfil#\qquad\cr
\multispan{12}\hfil\bf{Table 2}\hfil\cr
\multispan{12}\hfil Values of $P^{(1)}(n,I)$\hfil\cr
\noalign{\smallskip\hrule\smallskip}
n&I=0&1&2&3&4&5&6&7&8&9&10\cr
\noalign{\smallskip\hrule\smallskip}
0&0& & & & & & & & & & \cr
1&0& 0& & & & & & & & & \cr
2&0& 0& 0& & & & & & & & \cr
3&1& 1& 0& 0& & & & & & & \cr
4&6& 6& 2& 0& 0& & & & & & \cr
5&27& 27& 12& 3& 0& 0& & & & & \cr
6&110& 110& 55& 19& 4& 0& 0& & & & \cr
7&429& 429& 229& 91& 27& 5& 0& 0& & & \cr
8&1638& 1638& 912& 393& 136& 36& 6& 0& 0& & \cr
9&6188& 6188& 3549& 1614& 612& 191& 46& 7& 0& 0& \cr
10&23256& 23256& 13636& 6447& 2601& 897& 257& 57& 8& 0& 0\cr
\noalign{\smallskip\hrule}}}}
\bigskip

{}From this, one conjectures that $P^{(1)}(n,I)$ is given by the expression
$$g(n,I):=\binom{2n{-}I{-}1}n{-}\binom{2n{-}I{-}1}{n{+}3}{+}\binom{2n{-}2I{-}2}{n{-}I{-}4}{-}\binom{2n{-}2I{-}2}{n{-}I{-}1}{+}\binom{2n{-}2I{-}3}{n{-}I{-}4}{-}\binom{2n{-}2I{-}3}{n{-}I{-}2}.
$$
Since it is readily verified that $g(n,I)$ also satisfies the same recurrence and initial conditions, we have a rigorous proof that $P^{(1)}(n,I)=g(n,I)$. 
Plug in $I=1$ and we find that $a_{123}^{(1)}(n)=g(n,1)={3\over n}\binom{2n}{n+3}$
as first proved in \cite{4}.

\head 1.4 The number of permutations with exactly $2$ $abc$ patterns\endhead

We now will compute $a_{123}^{(2)}=$ the number of permutations in $S_n$ containing exactly two increasing subsequences of length 3. If 
$$
F_n^{123}=\sum_{\sigma\in S_n}q^{\varphi_{abc}(\sigma)}=a_{123}(n)+a_{123}^{(1)}(n)q+
a_{123}^{(2)}(n)q^2+\dots
$$
then 
$$
{{d^2F_n^{123}}\over{dq^2}}=
2a_{123}^{(2)}(n) +6a_{123}^{(3)}(n)q +\dots\ .
$$
So $a_{123}^{(2)}(n)$ is half the constant term of ${{d^2F_n^{123}}\over{dq^2}}$.
Similarly, by our definition of $P^{(r)}(n,I)$ (definition 1.5), we see that
$$
P^{(2)}(n,I)={1\over 2}\Phi^{(2)}(n,I).
$$
where 
$$
\displaystyle\Phi^{(2)}(n,I):={{\partial^2}\over {\partial q^2}}P_n\qqqline{q_k{\rightarrow }0,\  2{\leq}
k{<}i}{q_k{\rightarrow }1,\ i{\leq} 
k{\leq} n{-}1}{q{\rightarrow }0}.
$$
We will find a recursive formula for $\Phi^{(2)}(n,I)$. 

From (8), we have
$$
\eqalignno{
\quad\p2dir qP_n&(q_2,q_3,\dots,q_n;q)= 
{\partial^2\over{\partial
q^2}}P_{n{-}1}(q_3,q_4,\dots,q_n;q)&\cr&+\sum_{i=2}^n q_i^{i{-
}1}\left[ \p2dir qP_{n{-}
1}\qline{+}\sum_{j=2}^{i{-}1}q_j{\partial\over{\partial 
q_j}}\pdir qP_{n{-}1}\qline\right.&\cr&{+}\left.\sum_{j=2}^{i{-}1}q_j\left(\pdir q{\partial\over{\partial 
q_j}}P_{n{-}1}\qline{+}\sum_{m=2}^{i{-}1}q_m{\partial\over{\partial q_m}}{\partial\over{\partial q_j}}P_{n{-}1}\qline\right)\right].\quad&(13)}
$$
Let $\displaystyle\Phi^{(2)}_{(1,j)}(n,I)={\partial\over {\partial q}}\pdir {q_j}P_n\qqqline{q_k{\rightarrow }0,\  2{\leq}
k{<}i}{q_k{\rightarrow }1,\ i{\leq} 
k{\leq} n{-}1}{q{\rightarrow }0}
$ and $\displaystyle\Phi^{(2)}_{(j,m)}(n,I)={{\partial^2}\over {\partial q_j}}{\partial\over{\partial q_m}}P_n\qqqline{q_k{\rightarrow }0,\  2{\leq}
k{<}i}{q_k{\rightarrow }1,\ i{\leq} 
k{\leq} n{-}1}{q{\rightarrow }0}$. Then from (13) it follows that
$$
\Phi^{(2)}(n,I)=\sum_{i=I}^n\left[\Phi^{(2)}(n-1,i-1)+\sum_{j=I+1}^{i-1}\left(2\Phi^{(2)}_{(1,j)}(n-1,i-1)+\sum_{m=I+1}^{i-1}\Phi^{(2)}_{(j,m)}(n-1,i-1)\right)\right].
$$
Subtracting successive terms:
$$\eqalignno{
\Phi^{(2)}(n,I)-\Phi^{(2)}(n,I+1)&=\Phi^{(2)}(n-1,I-1)&\cr&+\sum_{i=I+2}^n\left[\Phi^{(2)}_{(I+1,I+1)}(n-1,i-1)+2\Phi^{(2)}_{(1,I+1)}(n-1,i-1)\right]&\cr&+\sum_{i=I+2}^n\sum_{j=I+2}^{i-1}2\Phi^{(2)}_{(I+1,j)}(n-1,i-1)&(14)\cr}
$$
We note now that to compute $\Phi^{(2)}(n,I)$ we must also compute $\Phi^{(2)}_{(1,j)}(n,I)$ and $\Phi^{(2)}_{(m,j)}(n,I)$, but we do not need to compute these for all $n,I,j,m$. We need only compute them for $j,m\leq I$. We use (8) to obtain a recursive formula for $\Phi^{(2)}_{(1,j)}(n,I)$, $j\leq I$.
$$
\eqalignno{
\quad\pdir q\pdir {q_j}P
_n&(q_2,q_3,\dots,q_n;q)&\cr&= \sum_{i=2}^{j-1}q_i^{i-1}\left[\sum_{m=2}^{i-1}q_m\pdir {q_m}\pdir{q_{j-1}}P_{n-1}\q1line i+\pdir q\pdir {q_{j-1}} P_{n-1}\q1line i\right] &\cr&+q\sum_{i=j+1}^n q_i^{i-1}\left[\sum_{m=2}^{i-1}q_m\pdir{q_m}\pdir{q_j}P_{n-1}\q1line i + \pdir q\pdir{q_j}P_{n-1}\q1line i\right]&\cr &+(j-1)q_j^{j-2}\left[\sum_{i=2}^{j-1}q_i\pdir{q_i}P_{n-1}\q1line j +\pdir q P_{n-1}\q1line j\right]&\cr&+\sum_{i=j+1}^nq_i^{i-1}\pdir {q_j}P_{n-1}\q1line i +\pdir q\pdir{q_{j-1}}P_{n-1}\qqline{q_k\rightarrow q_{k+1}}{q\rightarrow q}
&(15)}
$$
From this, we have
$$
\Phi^{(2)}_{(1,j)}(n,I)=\cases\sum_{i=I+1}^nP^{(1)}_{(2)}(n-1,i-1)+P^{(1)}(n-1,I-1),&\roman{if}\ j=2\cr
\sum_{i=I+1}^nP^{(1)}_{(j)}(n-1,i-1)+\Phi^{(2)}_{(1,j-1)}(n-1,I-1),&\roman{if}\ 2<j\leq I\cr\endcases
$$
where $P^{(1)}_{(j)}(n,I)=\pdir{q_j}P_n\q1line I$. This recurrence can be further simplified using (12) to
$$
\Phi^{(2)}_{(1,j)}=P^{(1)}(n-j+1,I-j+1)+\sum_{k=2}^j\left(\sum_{i=I+k-j+1}^{n-j+k}P(n-j,i-k)\right),\ \roman{for}\ j\leq I.\eqno(16)
$$
We may do the same for $\Phi^{(2)}_{(m,j)}(n,I)$ and obtain
$$
\Phi^{(2)}_{(m,j)}(n,I)=\cases
0,&\roman{if}\ j=m=2\cr
2P(n-1,I-1),&\roman{if}\ j=m=3\leq I\cr
\Phi^{(2)}_{(j-1,j-1)}(n-1,I-1),&\roman{if}\ j=m\leq I\cr
P^{(1)}_{(j-1)}(n-1,I-1),&\roman{if}\  j>m=2\cr
\Phi^{(2)}_{(m-1,j-1)}(n-1,I-1),&\roman{if}\  j>m,\ m\leq I\cr
\endcases
$$
which collapses down to the 3 cases
$$
\Phi^{(2)}_{(m,j)}(n,I)=\cases
0,&\roman{if}\ j=m=2\cr
2P(n-j+1,I-j+1),&\roman{if}\ j=m\leq I\cr
P^{(1)}_{(j-m+1)}(n-m+1,I-m+1),&\roman{if}\  j>m,\ m\leq I\cr
\endcases . \eqno(17)
$$
One might think that we have overlooked a few cases here, for what if $j>m>I$ or $j=m>I$? Examining (14) it is apparent that to compute $\Phi^{(2)}(n,I)$, we need only compute $\Phi^{(2)}_{(m,j)}$ and $\Phi^{(2)}_{(1,j)}$ when $j,m\leq I$. We now use (5),(12),(16) and (17) to simplify  (14):
$$
\eqalign{
\Phi^{(2)}(n,I)&=\Phi^{(2)}(n-1,I-1)+\Phi^{(2)}(n,I+1)\cr
&+2P^{(1)}(n-I,2)+2IP(n-I,3)+2P(n-i-1,2)+\chi(I>1)2P(n-I+1,3)\cr}
$$
We may now make the substitution $P^{(2)}(n,I)={1\over 2}\Phi^{(2)}(n,I)$. Our recurrence for $P^{(2)}(n,I)$ can be stated as
$$
P^{(2)}(n,I)=\cases
P^{(2)}(n,1),&\roman{if}\  I=0\cr
n-3,&\roman{if}\ I=n-2\cr
P^{(2)}(n-1,I-1)+P^{(2)}(n,I+1)+P^{(1)}(n-I,2)&\cr
\qquad+P(n-i-1,2)+IP(n-I,3)+\chi(I>1)P(n-I+1,3),&\roman{if}\ 0<I<n-2, n>3\cr
\endcases.
$$
\bigskip
\centerline{
\vbox{
\halign{#\hfil\qquad& \qquad\hfil#\qquad&\hfil#\qquad& \hfil#\qquad& 
\hfil#\qquad& \hfil
#\qquad& \hfil#\qquad& \hfil#\qquad& \hfil#\qquad& \hfil#\qquad& \hfil#\qquad& 
\hfil#\qquad\cr
\multispan{12}\hfil\bf{Table 3}\hfil\cr
\multispan{12}\hfil Values of $P^{(2)}(n,I)$\hfil\cr
\noalign{\smallskip\hrule\smallskip}
n&I=0&1&2&3&4&5&6&7&8&9&10\cr
\noalign{\smallskip\hrule\smallskip}
0&   0&  0&  0&  0& 0&0&0 & 0  &0 & 0  &0 \cr
1&   0&  0&  0&  0&  0&0&0 & 0  &0  &0  &0 \cr
2&   0&  0&  0&  0&  0& 0&0 & 0  &0  &0  &0 \cr
3&   0&  0&  0&  0&  0& 0&0& 0  &0  &0  &0 \cr
4&   3&  3&  1&  0&  0& 0&0&0  &0  &0  &0 \cr
5&   24& 24& 12& 2&  0& 0&0&0 & 0 &0  &0 \cr
6&  133&133& 74& 23& 3& 0&0&0 & 0  &0  & 0 \cr
7&  635&635&371&141&36& 4&0&0 & 0 & 0  &0  \cr
8&  2807 &  2807  & 1688   &709&227 &  51&5&0  &0 & 0 & 0 \cr
9& 11864 & 11864  & 7276  & 3248 & 1168 &  334  & 68 &  6  &0 & 0  &0 \cr
10& 48756 & 48756 & 30340 & 14121 & 5459  &1771 & 464  &87 & 7  &0 & 0   \cr
\noalign{\smallskip\hrule}}}}
\bigskip
Using {\tt ordi} of \cite{10} or {\tt gfun}\cite{5} we conjecture that 
$$
a_{123}^{(2)}(n)={{59n^2+117n+100}\over{2n(2n-1)(n+5)}}\binom{2n}{n-4}.
$$
It is very likely that one should be able to conjecture an explicit expression for 
$P^{(2)}(n,I)$, which would be routine to prove, and from which the above conjecture 
would follow.
\head 2 The Method\endhead
We may now outline the method described in this paper.\par
To determine the number of permutations on $\{1\dots n\}$ having exactly $r$ occurrences of the pattern $x_1x_2\dots x_k$,\par\noindent
\proclaim{1} Determine the best way to obtain a recurrence for this pattern. \endproclaim
There are basically four ways to do this.\par
\qquad a) by removing the last entry of the permutation, as in the $abc$ example.\par
\qquad b) by removing the first entry of the permutation, which is what we will do in the next example.\par
\qquad c) by removing $n$ from the permutation.\par
\qquad d) by removing $1$ from the permutation.\par
\noindent
\proclaim{2} Identify the other parameters needed in order to describe the recurrence.\endproclaim
In our first example, we found a recurrence for the number of permutations with a given number of $abc$ patterns by looking at what happened to a permutation from which we removed its last object. As a consequence, we were forced to consider the number of $abc$ patterns present in each permutation, but also the number of $aj$ patterns present. Note that we could arrive at this requirement by noting the result of removing the last object from the pattern $abc$. $abc$ becomes $ab$. Since only some choices of $b$ result in a true $abc$ pattern for a given $c$, we must be specific, and count the number of $aj$ patterns for every $j$. We count the number of $aj$ by using the parameter $q_j$ in our weight function. 
\proclaim{3} Define $P_n=\sum_{\sigma\in S_n}wt(\sigma)$.\endproclaim
Here $wt(\sigma)=q^\varphi(\sigma)\prod q_j^{\varphi_j(\sigma)}$ where $\varphi_j(\sigma)$ is the number of times the pattern associated with $q_j$ can be found in $\sigma$ and $\varphi(\sigma)$ is the number of $x_1x_2\dots x_r$ in $\sigma$.
\proclaim{4} Determine the functional equation.\endproclaim
Using the recurrence described in {\bf 1}, determine a recurrence in $n$ for $P_n$.
\proclaim{5} Take the $r^{th}$ derivative of the functional equation with respect to $q$.\endproclaim
\proclaim{6} Let $q=0$, and $q_j=1$ for all other parameters of $P_n$. \endproclaim
For all intents and purposes, you are now done, for you have obtained a (perhaps complicated) recurrence for the number of permutations on $\{1\dots n\}$ containing exactly $k$  $x_1x_2\dots x_r$'s. This recurrence probably involves other terms, but each of these have their own recurrences which can be determined from the functional equation. After all is said and done, you frequently are able to simplify this recurrence, eliminating many of these unwanted terms, but for now, you have achieved your goal.
\medskip
\head 3 The forbidden pattern $cab$\endhead
\head 3.1 Definitions\endhead
As with each of the examples we examine, the definitions for $P$, $P^{(k)}$, $\varphi_j$ and $wt$ should be taken as being local to the problem at hand.
\definition{Definition 3.1} Given $\sigma\in S_n$, a $cab$ pattern is a sequence $i,j,k$ where $1\leq i<j<k\leq n$ and $\sigma(j)<\sigma(k)<\sigma(i)$.\enddefinition
\definition{Definition 3.2}For $\sigma\in S_n$, let $\varphi_{cab}(\sigma)$ be the number of cab patterns of $\sigma$.\enddefinition
 The number of permutations on $\{1\dots n\}$ having no $cab$ patterns will be the constant term of the polynomial $F=\sum_{\sigma\in S_n}q^{\varphi_{cab}(\sigma)}$. As we saw earlier, we may add any number of parameters to this polynomial, as long as we know what to do with them to obtain the constant term. In this case, we will use the parameters $q_2,...,q_n$ as we did with $abc$. 
\definition{Definition 3.3}Given $\sigma\in S_n$, 
$$
wt(\sigma)=q^{\varphi_{cab}(\sigma)} \prod_{j=2}^nq^{\varphi_j(\sigma)}.
$$
\enddefinition
\head 3.2 No $cab$'s\endhead
Let $\sigma(n)=i$. A functional equation results from examining $\sigma_1(k)=\sigma(k),\ 1\leq k<n$,
$$
wt(\sigma)=wt(\sigma_1)q^{\sum_{j<i}\varphi_j(\sigma_1)}\prod_{j>i}q_j
$$
Summing over all $\sigma\in S_n$ we have
$$
P_n(q,q_2,q_3,\dots,q_n)=\sum_{i=2}^n\left[\left(\prod_{j>i}q_j\right)P_{n-1}(q,qq_2,\dots,qq_{i-1},q_{i+1},\dots,q_n)\right]+\left(\prod_jq_j\right)P_{n-1}(q,q_3,\dots,q_n).\eqno(18)
$$
Not surprisingly, if we let $P(n,I)=P_{n}(0;\overbrace{0,0,\dots ,0}^{I-1},\overbrace{1,\dots,1}^{n-I-1})$, then we obtain the following recurrence for $P(n,I)$,
$$
P(n,I) = \cases 1\hfill ,&\roman{if}\ n=I\hfill\cr
P(n,1)\hfill ,&\roman{if}\ I=0\hfill\cr
P(n,I+1)+P(n-1,I-1),&\roman{otherwise}\cr \endcases
$$
Which is identical to (6). So we have reproved that the number of permutations on $\{1\dots n\}$ with no $abc$ subsequences is equal to the number with no $cab$ 
subsequences.

\head 3.3 Permutations with one $cab$\endhead

Though the number of permutations with no $abc$'s is equal to the number with no $cab$'s, we will find that this is not the case when we examinee the number of permutations with one $cab$. we follow the same procedure outlined for the $abc$ case, with the goal of finding $\displaystyle \pdir q P_n\qqline{q\rightarrow 0}{q_j\rightarrow 1,\ 1\leq j\leq n}$.
Taking derivatives of both sides of $(18)$ with respect to $q$, we have
$$\eqalign{
\pdir qP_n&=\sum_{i=1}^n\left(\prod_{j>i}q_j\right)\left[\sum_{j=2}^{i-1}q_j\pdir{q_j}P_{n-1}\qline+\pdir{q}P_{n-1}\qline\right]\cr &+\left(\prod_jq_j\right)\pdir qP_{n-1}\qqline{q\rightarrow q}{q_k\rightarrow q_{k+1},\ 1\leq k\leq n-1}.\cr}
$$
We do the same for the derivative with respect to $q_j$:
$$\eqalign{
\pdir {q_j}P_n&=\sum_{i=2}^{j-1}\left(\prod_{{k>i}\atop{k\neq j}}q_k\right)P_{n-1}\qline+
\sum_{i=1}^{j-1}\left(\prod_{k>i}q_k\right)\pdir{q_{j-1}}P_{n-1}\qline \cr &+\sum_{i=j+1}^{n}\left(\prod_{k>i}q_k\right)q\pdir{q_{j}}P_{n-1}\qline.\cr}
$$
Here we let $P^{(1)}(n,I)=\pdir qP_n\qqqline{q\rightarrow 0}{q_k\rightarrow 0,\ k\leq I}{q_k\rightarrow 1,\ k>I}$ and $P^{(1)}_{(j)}(n,I)=\pdir {q_j}P_n\qqqline{q\rightarrow 0}{q_k\rightarrow 0,\ k\leq I}{q_k\rightarrow 1,\ k>I}$. We obtain the following recursions:
$$
\eqalign{
P^{(1)}_{(j)}(n,I)&=\cases 0,&\roman{if}\ n-I<0 \ or \ j<I\cr
P^{(1)}_{(j)}(n,1),&\roman{if}\ I=0\cr
P(n-1,I-1),&\roman{if}\ j=I\cr
\sum_{k=i}^{j-1}P(n-1,k-1),&\roman{if}\ j>I\cr\endcases ,\cr
&\cr
P^{(1)}(n,I)&=\cases I,&\roman{if}\ I=n-2\cr
P^{(1)}(n,1),&\roman{if}\ I=0\cr
\sum_{i=I}^n\left[\sum_{j=I+1}^{i-1}P^{(1)}_{(j)}(n-1,i-1)+P^{(1)}(n-1,i-1)\right],&\roman{otherwise}\cr\endcases .\cr}
$$
We combine these two recurrences and localize to obtain 
$$
P^{(1)}(n,I)=\cases I,&\roman{if}\ I=n-2\cr
P^{(1)}(n,1),&\roman{if}\ I=0\cr
P^{(1)}(n,I+1)+P^{(1)}(n-1,I-1)+P(n-2,I),&\roman{otherwise}\cr\endcases .
$$
\bigskip
\centerline{
\vbox{
\halign{#\hfil\qquad& \qquad\hfil#\qquad&\hfil#\qquad& \hfil#\qquad& 
\hfil#\qquad& \hfil
#\qquad& \hfil#\qquad& \hfil#\qquad& \hfil#\qquad& \hfil#\qquad& \hfil#\qquad& 
\hfil#\qquad\cr
\multispan{12}\hfil\bf{Table 4}\hfil\cr
\multispan{12}\hfil Values of $P^{(1)}(n,I)$\hfil\cr
\noalign{\smallskip\hrule\smallskip}
n&I=0&1&2&3&4&5&6&7&8&9&10\cr
\noalign{\smallskip\hrule\smallskip}
0&0& & & & & & & & & & \cr
1&0& 0& & & & & & & & & \cr
2&0& 0& 0& & & & & & & & \cr
 3&     1     & 1      &0  &    0     & &        & &     & &    & \cr
 4&     5      &5      &2   &   0     &0  &       &  &    & &   & \cr
 5&     21     &21     &11   &  3     &0   &  0    &   &   & &   & \cr
 6&     84     &84     &49    &19     &4    & 0    &0    &  & &   & \cr
 7&    330   & 330    &204    &92    &29     &5    &0    &0  & &   & \cr
 8&    1287   &1287   &825    &405   &153   &41    &6    &0  &0 &   & \cr
 9&    5005   &5005  & 3289  &1705  & 715   &235  & 55   &7  &0 & 0 &  \cr
 10&   19448 & 19448 & 13013  &7007 & 3146  &1166 & 341  &71 & 8 & 0 & 0 \cr
\noalign{\smallskip\hrule}}}}
\bigskip
Using {\tt ordi}\cite{10} or {\tt gfun}\cite{5} we conjecture that $a^{(1)}_{312}(n)$,
the number of permutations on $\{1,2,\dots,n\}$ containing exactly one $cab$ subsequence, $={{n-2}\over{2n}}\binom{2n-2}{n-1}$.
\head 4 Counting $abcd$'s\endhead
Lest the reader think that the methods outlined in this paper will only help us gain information about permutations avoiding forbidden patterns of length three, here we examine the forbidden pattern $abcd$. 
\head 4.1 Definitions\endhead
\definition{Definition 4.1}An $abcd$ pattern of a permutation $\sigma$ is a sequence $1\leq i<j<k<l\leq n$ with $\sigma(i)<\sigma(j)<\sigma(k)<\sigma(l)$. 
\enddefinition
\definition{Definition 4.2} $\varphi_{abcd}(\sigma)=$ the number of $abcd$ patterns which can be found in $\sigma$.\enddefinition
As before, we have
\definition{Definition 4.3}
$\varphi_{aj}(\sigma)=$ the number of $aj$ patterns of $\sigma$. \enddefinition
Here we must also define the following.
\definition{Definition 4.4} An $abj$ pattern is a sequence $1\leq i_1<i_2<i_3\leq n$ with $\sigma(i_1)<\sigma(i_2)<\sigma(i_3)=j$.\enddefinition
\definition{Definition 4.5}$\varphi_{abj}(\sigma)=$ the number of $abj$ patterns of $\sigma$.\enddefinition
 \definition{Definition 4.6}Let $\sigma\in S_n$ then 
$$
wt(\sigma)=q^{\varphi_{abcd}(\sigma)}\prod_{j=3}^nq_j^{\varphi_{abj}(\sigma)}\prod_{j=2}^n\xi_j^{\varphi_{aj}(\sigma)}. 
$$
\enddefinition
\head The number of permutations with no $abcd$\endhead
Let $\sigma_1(j)=\sigma(j),\ 1\leq j\leq n-1$. If $\sigma(n)=i$ then
$$
wt(\sigma)=wt(\sigma_1)q^{\sum_{j=3}^{i-1}\varphi_{abj}(\sigma_1)}q_i^{\sum_{j=2}^{i-1}\varphi_{aj}(\sigma_1)}\xi_i^{i-1}.\eqno(19)
$$
Let $P_{n}(q;q_3,\dots,q_n;\xi_2,\dots,\xi_n)=\sum_{\sigma\in S_n}wt(\sigma)$. Using $(19)$, we have
$$\eqalignno{
P_{n}(q;q_3,\dots,q_n;\xi_2,\dots,\xi_n)&=\sum_{i=3}^n\xi_i^{i-1}P_{n-1}(q;qq_3,\dots,qq_{i-1},q_{i+1},\dots,q_n;q_i\xi_2,\dots,q_i\xi_{i-1},\xi_{i+1},\dots,\xi_n)&\cr
&+\xi_2P_{n-1}(q;q_4,\dots,q_n;\xi_3,\dots,\xi_n)+P_{n-1}(q;q_4,\dots,q_n;\xi_3,\dots,\xi_n).&(20)\cr}
$$ 
The number of permutations on $\{1\dots n\}$ with no $abcd$ patterns is 
$$
P_n(0;1,\dots,1;1,\dots,1)=\sum_{i=3}^nP_{n-1}(0;\overbrace{0,\dots,0}^{i-2},\overbrace{1,\dots,1}^{n-i-1};1,\dots,1)+2P_{n-1}(0;1,\dots,1;1,\dots,1).
$$
We may use $(20)$ twice more to find 
$$\eqalign{
P_{n}(1;\overbrace{0,\dots,0}^{I_1-2},\overbrace{1,\dots,1}^{n-I_1-1};1,\dots,1)&=\sum_{i=3}^{I_1}P_{n-1}(0;\overbrace{0,\dots,0}^{I_1{-}3},\overbrace{1,\dots,1}^{n-I_1-1};\overbrace{0,\dots,0}^{i-2},\overbrace{1,\dots,1}^{n-i-1})\cr&+\sum_{i=I_1+1}^{n}P_{n-1}(0;\overbrace{0,\dots,0}^{i{-}3},\overbrace{1,\dots,1}^{n-i-1};1,\dots,1)\cr
&+2P_{n-1}(0;\overbrace{0,\dots,0}^{I_1-3},\overbrace{1,\dots,1}^{n-I_1-1};1,\dots,1)\cr}
$$
and
$$\eqalign{
P_{n}(1;\overbrace{0,\dots,0}^{I_1-2},\overbrace{1,\dots,1}^{n-I_1-1};\overbrace{0,\dots,0}^{I_2-1},\overbrace{1,\dots,1}^{n-I_2-1})&=\sum_{i=I_2+1}^{I_1}P_{n-1}(0;\overbrace{0,\dots,0}^{I_1{-}3},\overbrace{1,\dots,1}^{n-I_1-1};\overbrace{0,\dots,0}^{i-2},\overbrace{1,\dots,1}^{n-i-1})\cr&+\sum_{i=I_1}^{n}P_{n-1}(0;\overbrace{0,\dots,0}^{i{-}3},\overbrace{1,\dots,1}^{n-i-1};\overbrace{0,\dots,0}^{I_2-2},\overbrace{1,\dots,1}^{n-I_2-1})\cr
&+P_{n-1}(0;\overbrace{0,\dots,0}^{I_1-3},\overbrace{1,\dots,1}^{n-I_1-1};\overbrace{0,\dots,0}^{I_2-2},\overbrace{1,\dots,1}^{n-I_2-1}).\cr}\eqno(21)
$$
Let $P(n,I_1,I_2)=P_{n}(1;\overbrace{0,\dots,0}^{I_1-2},\overbrace{1,\dots,1}^{n-I_1-1};\overbrace{0,\dots,0}^{I_2-1},\overbrace{1,\dots,1}^{n-I_2-1})$ for $0<\{I_1,I_2\}\leq n$. Observe that $P(n,n,n)$ is the number of permutations on $\{1\dots n\}$ with no $abcd$ pattern, no $abc$ pattern (for any $c$) and no $aj$ patterns for any $j$. There is only one such permutation, namely $[n,n-1,\dots,2,1]$, thus $P(n,n,n)=1$ for all $n$. Furthermore, it is clear from $(21)$ that when $I_1<I_2$, $P(n,I_1,I_2)=P(n,I_2,I_2)$. We may define $P(n,I_1,0)=P(n,I_1,1)$ and $P(n,0,I2)=P(n,1,I2)$ for all values of $n,\ I_1$ and $I_2$. Using this notation, we have
$$
P(n,I_1,I_2)=\cases P(n,1,I_2),&\roman{if}\ I_1=0\cr
P(n,I_1,1),&\roman{if}\ I_2=0\cr
P(n,I_2,I_2),&\roman{if}\ I_1<I_2\cr
1,&\roman{if}\ I_1=I_2=n\cr
\eqalign{\displaystyle{\sum_{i=I_2+1}^{I_1}}P(n-1,I_1-1,i-1)&+\displaystyle{\sum_{i=I_1+1}^n}P(n-1,i-1,I_2)\cr&+P(n-1,I_1-1,I_2-1)\cr},&\roman{otherwise}\endcases.
$$
\bigskip
\centerline{
\vbox{
\halign{\hfil#\qquad& \hfil#\quad&\hfil#\quad& \hfil#\quad& 
\hfil#\quad& \hfil
#\quad& \hfil#\quad& \hfil#\quad& \hfil#\quad& \hfil#\quad& \hfil#\quad& 
\hfil#\quad\cr
\multispan{12}\hfil\bf{Table 5}\hfil\cr
\multispan{12}\hfil Values of $P(n,I,1)$\hfil\cr
\noalign{\smallskip\hrule\smallskip}
n&I=0&1&2&3&4&5&6&7&8&9&10\cr
\noalign{\smallskip\hrule\smallskip}
0&0& & & & & & & & & & \cr
1&1& 1& & & & & & & & & \cr
2&2& 2& 2& & & & & & & & \cr
3&6& 6& 6& 5& & & & & & & \cr
4&23& 23& 23& 20& 14& & & & & & \cr
5&103& 103& 103& 92& 70& 42& & & & & \cr
6&513& 513& 513& 466& 372& 252& 132& & & & \cr
7&2761& 2761& 2761& 2536& 2086& 1509& 924& 429& & & \cr
8&15767& 15767& 15767& 14594& 12248& 9227& 6127& 3432& 1430& & \cr
9&94359& 94359& 94359& 87830& 74772& 57894& 40403& 24882& 12870& 4862& \cr
10&586590& 586590& 586590& 548325& 471795& 372565& 268909& 175474& 101036& 48620& 16796\cr
\noalign{\smallskip\hrule}}}}
\bigskip

\bigskip
\centerline{
\vbox{
\halign{\hfil#\qquad& \quad\hfil#\quad&\hfil#\quad& \hfil#\quad& 
\hfil#\quad& \hfil
#\quad& \hfil#\quad& \hfil#\quad& \hfil#\quad& \hfil#\quad& \hfil#\quad& 
\hfil#\quad\cr
\multispan{12}\hfil\bf{Table 6}\hfil\cr
\multispan{12}\hfil Values of $P(n,1,I)$\hfil\cr
\noalign{\smallskip\hrule\smallskip}
n&I=0&1&2&3&4&5&6&7&8&9&10\cr
\noalign{\smallskip\hrule\smallskip}
0&1& & & & & & & & & & \cr
1&1&1& & & & & & & & & \cr
2&2&2&1& & & & & & & & \cr
3&6&6&3&1& & & & & & & \cr
4&23&23&12&4&1& & & & & & \cr
5&103&103&56&20&5&1& & & & & \cr
6&513&513&288&110&30&6&1& & & & \cr
7&2761&2761&1588&640&190&42&7&1& & & \cr
8&15767&15767&9238&3882&1235&301&56&8&1& & \cr
9&94359&94359&56094&24358&8187&2163&448&72&9&1& \cr
10&586590&586590&352795&157265&55235&15575&3528&636&90&10&1\cr
\noalign{\smallskip\hrule}}}}
\bigskip
Using {\tt ordi}\cite{10}, we conjecture that $a^{(0)}_{1234}(n)=P(n,1,1)$ satisfies the following recurrence with \penalty-10000$a^{(0)}_{1234}(0)=1$ and $a^{(0)}_{1234}(1)=1$:
$$
a^{(0)}_{1234}(n+2)=-9{{(n+1)^2}\over{(n+4)^2}}a^{(0)}_{1234}(n)+{{10n^2+42n+41}\over{(n+4)^2}}a^{(0)}_{1234}(n+1)
$$

\head 5 Counting permutations avoiding more than one forbidden pattern\endhead
The method outlined in this paper can be used to find recurrences for the number of permutations avoiding 2 or more forbidden patterns. The method is essentially the same, though you have to keep track of more parameters, and recurrences can be more complicated. In the example that follows, we seek the number of permutations on $\{1\dots n\}$ avoiding both $abc$ and $bac$. Let $\varphi_{bac}(\sigma)$ be the number of $bac$ patterns of the permutation $\sigma$. Let $\varphi_{ja}(\sigma)$ be the number of $ja$ patterns of $\sigma$ (that is the number of inversions `caused' by $j$). Let
$$
wt(\sigma)=q^{\varphi_{abc}(\sigma)}\xi^{\varphi_{bac}(\sigma)}\prod_{j=2}^n\left(q_j^{\varphi_{aj}(\sigma)}\xi_j^{\varphi_{ja}(\sigma)}\right).
$$
It may seem a bit wasteful to spend time defining both $\varphi_{aj}(\sigma)$ and $\varphi_{ja}(\sigma)$ when it is clear that $\varphi_{aj}(\sigma)+\varphi_{ja}(\sigma)=n-1$ but we will see that this `complication' in addition to the introduction of the parameters $\xi_j$ will pay off in the end.\par
$$\eqalign{
P_n(q;\xi;q_2,\dots,q_n;\xi_2,\dots,\xi_n)&=\sum_{\sigma\in S_n}wt(\sigma)=\sum_{i=1}^n\sum_{{\sigma\in S_n}\atop{\sigma(n)=i}}wt(\sigma)\cr
=\sum_{i=2}^nq_i^{i-1}\left(\prod_{j>i}\xi_j\right)&P_{n-1}(q;\xi;qq_2,\dots,qq_{i-1},q_{i+1},\dots,q_n;\xi\xi_2,\dots,\xi\xi_{i-1},\xi_{i+1},\dots,\xi_n)\cr
&+\left(\prod_{j}\xi_j\right)P_{n-1}(q;\xi;q_3,\dots,q_n;\xi_3,\dots,\xi_n).}
$$
Let $P(n,I)=P_n(0;0;\overbrace{0,\dots,0}^{I-1},\overbrace{1,\dots,1}^{n-I};\overbrace{0,\dots,0}^{I-1},\overbrace{1,\dots,1}^{n-I})$. Note that from our definition, $P(n,I)=0$ when $I>1$. Indeed, $P(n,I)$ is the number of permutations on $\{1\dots n\}$ containing no $abc$, no $bac$, no $aj$ for $j\leq I$ and no $ja$ for $j\leq I$. If there was a permutation for which this was true for an I>1 then we merely examinee the positions of the object 1 and 2 in the permutation. Let $\sigma(i)=1$ and $\sigma(j)=2$. If $i<j$ then $\sigma$ has an $a2$ pattern. If $j<i$ then $\sigma$ has a $2a$ pattern. Thus, $\sigma$ does not meet the requirements and so $P(n,I)=0$ for $I>1$. We have the recurrence
$$
P(n,I)=\cases P(n,1),&\roman{if}\ I=0\cr
0,&\roman{if}\ I>1\cr
1,&\roman{if}\ n=I=1\cr
\sum_{i=I}^nP(n-1,i-1),& \roman{otherwise}\cr\endcases.
$$
This recurrence reduces to 
$$
P(n,1)=\cases 1,&\roman{if}\ n=1\cr
2P(n-1,1),&\roman{if}\ n>1\cr\endcases.
$$
So we have reproved the well know result that the number of permutations on $\{1\dots n\}$ with no $abc$ and no $bac$ is $2^n$.\par

\head Conclusion \endhead

We saw in the above examples that in order to compute the quantity of interest, $a_n$, say,
we naturally introduced extra parameters $I_1, I_2, \dots$, and
a new quantity $F(n,I_1,I_2, \dots )$ such that $a_n$ was the special case
$F(n,1,1, \dots)$, in which all  the extra parameters $I_1,I_2, \dots$,
are set to $1$. Since the system of linear recurrence equations always
seems to have constant coefficients, the generating function in all
the corresponding continuous variables

$$
\tilde F(z;x_1,x_2,\dots):=\sum_{n,I_1,I_2, \dots } F(n,I_1,I_2,\dots) z^n x_1^{I_1} x_2^{I_2} \dots
$$
should be holonomic (multi-D-finite) in all its variables (because of the
nonstandard boundary conditions, it is not always a rational function).
It follows from the holonomic
theory \cite{9} that any coefficient with respect to $x_1,x_2,\dots$, in particular that of $x_1^1x_2^1\dots$, is still holonomic (D-finite, i.e. satisfies a linear differential equation with polynomial coefficients), hence $f_n:=F(n,1,1,\dots)$ is $P-$recursive.

The method described here works for many patterns and sets of patterns, but it does not seem to work for all patterns. The authors were unable to find a suitable set of parameters (see the method, step 2) to apply this method to the pattern $adcb$.

\head Appendix\endhead
A maple package which confirms and illustrates many results from this paper is available and can be obtained using your favorite world wide web browser at {\tt http://www.math.temple.edu/\~{}noonan} or \penalty-10000{\tt http://www.math.temple.edu/\~{}zeilberg}.

\enddocument

\bye